\documentclass[12pt]{article}
\usepackage{amsmath,amsthm,amsfonts,latexsym,amssymb,enumerate,color}
\usepackage{graphicx}

\def\CC{\mathbb C}
\def\DD{\mathbb D}
\def\M{\mathcal M}
\def\RR{\mathbb R}
\def\TT{\mathbb T}

\def\H2p{\overline{H^2_0}} 
\def\Hlp{H^2(\CC^-)}
\def\cprime{$'$}

\def\beq{\begin{equation}}
\def\eeq{\end{equation}}
\def\ds{\displaystyle}

\def\Hol{\mathop{\rm Hol}\nolimits}
\def\ker{\mathop{\rm ker}\nolimits}
\def\Kmin{\mathop{\rm K_{min}}\nolimits}

\def\spam{\mathop{\rm span}\nolimits} 

\newtheorem{thm}{Theorem}[section]
\newtheorem{prop}[thm]{Proposition}
\newtheorem{defn}[thm]{Definition}
\newtheorem{lem}[thm]{Lemma}
\newtheorem{cor}[thm]{Corollary}
\newtheorem{rem}[thm]{Remark}

\newtheorem{ex}[thm]{Example}

\numberwithin{equation}{section}

\def\beginpf{\begin{proof}}
\def\endpf{\end{proof}}

\begin{document}

\title{Multipliers and equivalences between Toeplitz kernels}


\author{M.~Cristina C\^amara\thanks{
Center for Mathematical Analysis, Geometry and Dynamical Systems,
Instituto Superior T\'ecnico, Universidade de Lisboa, 
Av. Rovisco Pais, 1049-001 Lisboa, Portugal.
 \tt ccamara@math.ist.utl.pt} \ and  Jonathan R.~Partington\thanks{School of Mathematics, University of Leeds, Leeds LS2~9JT, U.K. {\tt j.r.partington@leeds.ac.uk}}\ \thanks{Corresponding author}
}

\maketitle

\begin{abstract}
Multipliers between kernels of Toeplitz operators are characterised in terms of test functions (so-called maximal vectors for
the kernels);  these maximal vectors
may easily be parametrised in terms of inner and outer factorizations. Immediate applications
to model spaces are derived. The case of surjective multipliers is also analysed. These ideas are applied to describing
equivalences between two Toeplitz kernels.
\end{abstract}



\maketitle

\section{Introduction}

The starting point for this work is a result of Fricain, Hartmann and Ross \cite{FHR}, 
which gives a necessary and sufficient condition for a function $g$ to multiply a model space 
$K_\theta$ into another model space $K_\phi$ (all notation and definitions will be given later in this section).
This in turn was motivated by a more restrictive version of this question due to
Crofoot \cite{crofoot}.

The main result of \cite{FHR} says that $w$ multiplies $K_\theta$ into $K_\phi$ if and only if:\\
(i) $w$ multiplies the function $S^*\theta$ into $K_\phi$ (here $S^*$ denotes the backward shift), and\\
(ii) $w$ multiplies $K_\theta$ into $H^2$ (this may be expressed as a Carleson measure condition).\\

Now model spaces are kernels of particular Toeplitz operators, indeed $K_\theta= \ker T_{\overline\theta}$,
and thus the question may be posed more generally for kernels of Toeplitz operators. We may
also ask whether more general test functions can be used, other than $S^*\theta$.

In this paper we address these questions, obtaining the result above as an immediate corollary. 
To do this we need to bring in some of the theory of Toeplitz kernels, particularly ideas
developed by the authors in \cite{CP14,CMP}. That work was done in the context of
Hardy spaces on the half-plane, and we reformulate it for the disc, showing also how
 the multiplier problem is solved for the half-plane.

In Section \ref{sec:2}, we establish the notion of minimal kernels and maximal vectors for kernels of Toeplitz
operators on $H^2$, and then use these to give a characterization of multipliers from one
Toeplitz kernel to another by using the maximal vectors as test functions. From this
we easily recover results on model spaces as special cases.

We also use the theory of multipliers to obtain results on the structure of Toeplitz kernels, linked to
factorization results for their symbols, together with theorems linking
an equivalence between kernels with an
equvalence between their symbols. 

In Section \ref{sec:3}, we obtain necessary and sufficient conditions for {\em surjective\/} multipliers between Toeplitz kernels,
recovering Crofoot's result as a very special case.

In Section \ref{sec:4}, we give a brief discussion of the situation for the upper half-plane, which can
be obtained independently or by using the unitary equivalence of the corresponding Hardy spaces.

\subsection*{Notation}

We use $H^2$ to denote the standard Hardy space of the unit disc $\DD$, which embeds isometrically
into $L^2(\TT)$, where $\TT$ denotes the unit circle with normalized Lebesgue measure $m$. Its orthogonal complement is 
written $\overline{H^2_0}$ or $\overline{zH^2}$. Here $z$ denotes the independent variable.
The space $H^\infty$ is the Banach algebra of
bounded analytic functions on $\DD$, of which the set of invertible elements will be denoted by
$\mathcal G H^\infty$.
Moreover, $\Hol(\DD)$ denotes the space of all analytic functions on $\DD$.

We refer the reader to \cite{duren,hoffman,koosis,nik} for standard results on
Hardy spaces and the factorization of Hardy-class functions into inner and outer factors.

An observation that we shall use several times is that $f \in H^2$ if and only if $\overline z \overline f \in \overline{H^2_0}$, and likewise
$f \in \overline{H^2_0}$ if and only if $\overline z \overline f \in H^2$. 

The shift operator $S: H^2 \to H^2$ is the operator of multiplication by the independent variable $z$.

The Toeplitz operator $T_g$ with symbol $g \in L^\infty(\TT)$ is the operator on $H^2$ defined by
$T_g f= P_{H^2} (g   f)$, for $f \in H^2$, where $P_{H^2}$ denotes the orthogonal projection from
$L^2(\TT)$ onto $H^2$. 
If $\theta$ is an inner function, then $\ker T_{\overline\theta}$ is the model space $K_\theta= H^2 \ominus \theta H^2= H^2 \cap \theta \overline{H^2_0}$,
which is invariant under the backward shift $S^*$.

For $g, h \in L^\infty=L^\infty(\TT)$ we write 
 $\M(\ker T_g, \ker T_h)$ for the space of multipliers $w \in \Hol(\DD)$
such that $wf \in \ker T_h$ for all $f \in \ker T_g$ and we use the notation $\M_\infty(\ker T_g, \ker T_h)=\M(\ker T_g, \ker T_h)\cap L^\infty(\TT)$ and $\M_2(\ker T_g, \ker T_h)=\M(\ker T_g, \ker T_h)\cap L^2(\TT)$. \\

In fact, as we shall see later (Remark \ref{rem:24rev}), the multipliers between model spaces are necessarily contained in $H^2$;
this is not the case for general Toeplitz kernels, although they must lie in the Smirnov class.


\section{Multipliers and maximal vectors}
\label{sec:2}

\begin{defn}
For a function $k \in H^2 \setminus \{0\}$ we write $\Kmin(k)$ for the minimal Toeplitz kernel
containing $k$; that is, $\Kmin(k)= \ker T_v$ for some $v \in L^\infty $, with $k \in \Kmin(k)$, while 
 $\ker T_v \subset \ker T_w$ for
every $w \in L^\infty$ such that $k \in \ker T_w$.

We say that $k$ is a {\em maximal vector\/} for $\ker T_g$ if $\ker T_g = \Kmin(k)$.
\end{defn}

The existence of minimal kernels and maximal vectors was established in \cite[Thm 5.1 and Cor 5.1]{CP14}
in the context of the upper half-plane. Let us sketch the corresponding argument for the disc.\\

Suppose that $k=\theta p$, where $\theta$ is inner and $p$ is outer. Then we assert that
$\Kmin(k)=\ker T_v$, where $v=\overline z \overline \theta \overline p/p$.
Since $vk=\overline z \overline p$, we have $k \in \ker T_v$.

Now suppose that $k \in \ker T_w$ for some $w \in L^\infty $, and that $g \in \ker T_v$.
Thus $gv \in \overline{H^2_0}$ and $k w \in \overline{H^2_0}$.

Then $gw = gvkw/(vk)= (gv)(kw)/(\overline z \overline p)$; that is, $gw$ lies in   $L^2$,
and $\overline{zgw}=\overline{zgv}\overline{zkw}/p$, which means that $\overline{zgw}$ is in the
Smirnov class 
(the ratio of an $H^1$ function and an outer $H^2$ function) as well as $L^2(\TT)$.
By the generalized maximum principle (e.g. \cite[Thm. 2.11]{duren},\cite[Thm.~4.4.5]{nik}) 
it is therefore
in $H^2$ . Thus $gw \in \overline{H^2_0}$ and
$g \in \ker T_w$, and so $\Kmin(k)=\ker T_v$.

Moreover, by \cite[Lemma 1]{sarason94}, every Toeplitz kernel $K$ is $\ker T_{\overline z \overline \theta \overline p/p}$
for some inner function $\theta$ and outer function $p$ and thus $K=\Kmin(\theta p)$.\\

In fact, we can characterise all the maximal vectors for a Toeplitz kernel,  as follows.

\begin{thm}\label{thm:maxtg}
Let $g \in L^\infty\setminus\{0\}$ be such that $\ker T_g$ is non-trivial. Then
$k$ is a maximal vector for $\ker T_g$ if and only if $k \in H^2$ and
$k=g^{-1} \overline z \,\overline p$, where $p$
is outer in $H^2$.
\end{thm}

\beginpf
Note first that if $\ker T_g$ is non-trivial, then $gf \in \H2p$ for some nonzero $f \in H^2$, and so $g \ne 0$
almost everywhere and we can define $g^{-1}$.

Now if $\Kmin(k)=\ker T_g$, then we have $gk = \overline z \overline p$, where $p \in H^2$. 
Also $p$ is outer,
since if  $p=\phi q$, where $\phi$ is inner and non-constant, and $q$ is outer, then
 $k \in \ker T_{  \phi g} \subsetneq \ker T_g$, which contradicts the assumption.\\

Conversely, if $k=g^{-1} \overline z \,\overline p$, where $p$
is outer, then $k \in \ker T_g$. If also $k \in \ker T_h$ with $h \in L^\infty$, then
$\overline z \overline{hk} \in H^2$, and if $f \in \ker T_g$ we have $g f \in \H2p$, so
$\overline z \overline g \overline f \in H^2$.

Then 
\[
\overline z \overline{hf} = \overline z \overline{hk} \frac{\overline f}{\overline k}
= \overline z \overline{hk} \frac{\overline z \overline g \overline f}{\overline z\overline g \overline k}
=\overline z \overline{hk} \frac{\overline z \overline g \overline f}{p},
\]
which is in $L^2(\TT)$ and the Smirnov class, hence in $H^2$.
Thus $hf \in \H2p$ and $f \in \ker  T_h$; so $\ker T_g \subset \ker T_h$ and $\ker T_g = \Kmin(k)$.
\endpf

In the special case of a model space, we obtain immediately a
disc version of \cite[Thm. 5.2]{CMP}.
 
\begin{cor}
Let $\theta$ be inner. Then
$K_\theta=\Kmin(k)$ if and only if $k \in H^2$ and $k=\theta \overline z \overline p$, where $p$
is outer in $H^2$.
\end{cor}
\beginpf
Take $g=\overline\theta$ and apply Theorem \ref{thm:maxtg}.
\endpf

We are now ready to state a theorem characterizing multipliers of Toeplitz kernels.
Recall that $\mu$ is a Carleson measure for a subspace $X$ of $H^2$ if
there is a constant $C>0$ such that
\[
\int_{\TT} |f|^2 \, d\mu \le C \|f\|^2_2 \qquad \hbox{for all} \quad f \in X.
\]
In fact the measures that arise here will be supported on $\TT$, not $\DD$, and be absolutely continuous
with respect to Lebesgue measure, but it is convenient to see them in this more general perspective.
The natural choices for $X$ will be Toeplitz kernels, including model spaces. \\

Carleson measures for $\ker T_g$ may  be better understood if we use the
fact that $\ker T_g$ is nearly invariant, and thus 
by Hitt's result \cite{hitt}
$\ker T_g=FK_\theta$ for some isometric multiplier $F$
(which is outer) and $\theta$ inner.

We require $w$ to satisfy
\[
\|w F k\|_2 \le C \|Fk\|_2 = C \|k\|_2 
\] for each $k \in K_\theta$.
Thus the study of Carleson measures for Toeplitz kernels reduces to that of the special case
where the Toeplitz kernel is a
model space.
There is information on how to find an appropriate $\theta$ in Sarason's paper \cite{sarason94}.

Descriptions of Carleson measures for certain model spaces were given in \cite{cohn,VT},
with a complete answer in a recent preprint \cite{LSSUTW}.

We say that $w\in \mathcal C(\ker T_v)$ whenever $|w^2|dm$ is a Carleson measure for $\ker T_g$, that is $w\,\ker T_g \subset L^2(\TT)$.

\begin{rem}\label{rem:24rev}
{\rm
Note that   every 
nontrivial Toeplitz kernel contains an 
outer function, because if $\theta p \in \ker T_g$, where $\theta$ is inner and $p$ is
outer, then $p \in \ker T_g$ since $gp=\overline\theta (g\theta p) \in \overline{H^2_0}$.
Hence multipliers must   be holomorphic in $\DD$, and indeed lie in the Smirnov class $\mathcal N_+$. 
Moreover, a multiplier $w$ from a model space $K_\theta$, where $\theta$ is an inner function, into another Toeplitz kernel must be in $H^2$, 
since we must have $w\,(1-\overline{\theta (0)}\,\theta)\in H^2$, and $1-\overline{\theta (0)}\,\theta$ is invertible in $H^\infty$. 

Since Toeplitz kernels have the near-invariance property that $\theta p \in \ker T_g$ implies that
$p \in \ker T_g$, it follows easily that the space of multipliers has a similar property. Thus a non-zero multiplier space contains an outer function.

However, note that multipliers  between two general Toeplitz kernels
need not lie in $H^2. $ 
For example,
the function $z \mapsto (z-1)^{1/2}$ spans a 1-dimensional Toeplitz kernel $\ker T_g$, where $g(z)=z^{-3/2}$ with
$\arg z \in [0,2\pi)$ on $\TT$. This can be shown directly, or by using known results on the 
half-plane from \cite{CP14} together with the methods of Section \ref{sec:4} below. Hence the 
function $w(z)=(z-1)^{-1/2}$ multiplies $\ker T_g$ onto  the model space $K_z=\ker T_{\bar z}$ consisting only of the
constant functions, although $w$ is not an $H^2$ function. {It is easy to see that in fact $w$ satisfies conditions (ii) and (iii) in the following theorem.
}}
\end{rem}

\begin{thm}\label{thm:main}
Let $g, h \in L^\infty(\TT) \setminus \{0\}$ such that $\ker T_g$ and $\ker T_h$ are nontrivial.  
Then 
the following are equivalent:\\
(i)\,$w \in  \M(\ker T_g, \ker T_h);$\\
(ii)\,$w \in  \mathcal C(\ker T_g)$ and $w k \in \ker T_h$ for some (and hence all) maximal vectors $k$ of $\ker T_g$;\\
(iii)\,$w \in  \mathcal C(\ker T_g)$ and $h g^{-1}w\in \overline {\mathcal N_+}$.

\end{thm}

\beginpf
First we prove that (i)$\Leftrightarrow$(ii).
Clearly, the two conditions in (ii) are necessary  for (i). So assume that 
(ii) holds, and write $k=\theta p$, where
$\theta$ is inner and $p$ is outer.
Now $\ker T_g= \ker T_{\overline z\overline \theta\overline p/p}$, as detailed above,
and
thus without loss of generality we may take $g=\overline z\overline \theta\overline p/p$.

We have that $wkh \in \H2p$, since $wk \in \ker T_h$. Suppose now that
$f\in \ker T_g$, so that $ f g \in\H2p$.
Now 
\[
wf  h = (wkh) \frac{f}{\theta p} = (wkh) \frac{zfg}{\overline p}.
\]
Then $w f h \in L^2(\TT)$, since $wf \in L^2(\TT)$ by the Carleson condition.
Also $wkh$ and $f g$ are in $\H2p$ so  
$\overline{z}\overline{wfh}=\overline z \overline{wkh}\,\overline z \overline{fg}/p$
is in the Smirnov class
of the disc as well as $L^2(\TT)$.
Once again, we deduce that
$\overline{z}\overline{wfh} \in H^2$ and so $wf h \in \H2p$,
and finally $wf \in \ker T_h$.\\
Let now $w \in  \mathcal C(\ker T_g)$. To show that (ii)$\Rightarrow$(iii), assume that $k$ is a maximal vector for $\ker T_g$; then by Theorem \ref{thm:maxtg} we have $k=g^{-1}\bar z\bar p$ where $p$ is outer in $H^2$. If $w\ker T_g \subset \ker T_h$, then 
\[
hwk=hwg^{-1}\bar z\bar p=\psi_- \in \overline{H^2_0}
\]
so $hwg^{-1}=z \frac{\psi_-}{\bar p}\in \overline{\mathcal N_+}$.\\
Conversely, if $hwg^{-1}\in \overline{\mathcal N_+}$ then, for any maximal function $k$ of $\ker T_g$, for which  $gk \in \overline{H^2_0}$, we have
\[
h(wk)=hwg^{-1}(gk) \in \bar z \overline{\mathcal N_+}\cap L^2(\TT)=\overline{H^2_0}
\]
so $wk\in \ker T_h.$ 
\endpf

When $g=h$ and $\bar g$ is an inner function $\theta$, from Theorem \ref{thm:main} we get the well-known result that $\mathcal M (K_\theta, K_\theta) = \mathbb C$.

Note that if $k$ is not a maximal vector of $\ker T_g$, then $k$ cannot  be used as a test
function for multipliers from $\ker T_g$; for example in this case the   function $w(z) \equiv 1$ is not a multiplier from
$\ker T_g$ into $\Kmin(k)$, even though $wk \in \Kmin(k)$.

\begin{cor}\label{2.3A}
With the same assumptions as in Theorem \ref {thm:main}, and assuming moreover that $hg^{-1}\in L^\infty (\TT)$,
\[
w\in \M_2 (\ker T_g,T_h) \Leftrightarrow w \in  \mathcal C(\ker T_g) \cap \ker T_{\bar z h g^{-1}}.
\]
\end{cor}
\beginpf
Assume that $w\in \M_2 (\ker T_g,T_h)$; then $w\in H^2$ and from Theorem \ref {thm:main}(iii) it follows that $w \in  \mathcal C(\ker T_g)$ and $\bar z hg^{-1}w \in \overline{H^2_0}$, so that $w \in \ker T_{\bar z h g^{-1}}$. Conversely, if $w \in \ker T_{\bar z h g^{-1}}$ then $hg^{-1}w \in 
\overline{H^2} \subset \overline{\mathcal N_+}$, and the result follows from Theorem \ref{thm:main}.
\endpf
Regarding the assumption that $hg^{-1}\in L^\infty (\TT)$ in the corollary above, note that by [\cite {sarason94}, Lemma 1], for every Toeplitz kernel $K$ there exists $g \in L^\infty(\TT)$ with $|g|=1$ a.e. such that $K=\ker T_g$.

By considering in particular $g=\bar \theta$, where $\theta$ is an inner function, we obtain the following, which slightly generalises a result in \cite{FHR}.

\begin{cor} \label{cor:fhr-1}
Let $\theta$ be inner and $h\in L^\infty(\TT)\setminus\{0\}$ such that $\ker T_h$ is nontrivial. Then the following are equivalent:\\
(i)\,  $w \in \M(K_\theta,\ker T_h);$\\
(ii)\,  $wS^*\theta \in \ker T_h$, and $w \in  \mathcal C(K_\theta);$\\
(iii)\, $w \in  \ker T_{\bar z\theta h}\cap\mathcal C(K_\theta).$
\end{cor}

\beginpf
Since $S^*\theta = \theta \overline z \overline p$, where $p=1-\overline{\theta(0) }\theta$,
which is outer, we see
that $K_\theta=\Kmin(S^*\theta)$.
Thus the equivalence of (i) and (ii) follows directly from
Theorem \ref{thm:main}.

Finally, note that the first condition in (ii) asserts that $h w S^*\theta \in \H2p$ and  $w \in  \ker T_{\bar z\theta h}$ asserts that $ h w\theta\overline z
\in \H2p$. These conditions are  equivalent since $S^*\theta = \theta \overline z\overline{ (1-\overline{\theta (0)}\theta)}$,
where the last factor is invertible in $\overline{H^\infty}$.
\endpf




Note that, unlike $S^*\theta$, the reproducing kernel used as a test function in many other contexts,
beginning perhaps with \cite{bonsall},
is not maximal for $K_\theta$. For with
\[
k_a(z)=\frac{1-\overline{\theta(a)}\theta (z)}{1-\overline a z},
\]
we have 
\[
\theta \overline z \overline{k_a(z)} = \frac{\theta(z)-\ {\theta(a)} }{z- a  },
\]
which is not outer in general.

Corollaries \ref{2.3A} and \ref{cor:fhr-1} bring out a close connection between the existence of non-zero multipliers in $L^2(\TT)$ and their description, on the one hand, and the question of  injectivity of an associated Toeplitz operator $T_{\bar z\,g^{-1} h}$ (or $T_{\bar z\theta h}$) and the characterisation of its kernel, on the other hand.

It is well known that various properties of Toeplitz operators, in particular Toeplitz kernels, can be described in terms of a factorisation of their symbols. 

Recall that a function $f \in H^p \setminus \{0\}$ with $0<p<\infty$ is said to be {\em rigid}, if for any
$g \in H^p$ with $g/f > 0$ on $\TT$ we have $g=\lambda f$ for some $\lambda>0$. A rigid function is
outer, and every rigid function in $H^p$ is the square of an outer function in $H^{2p}$. A function $f \in H^2$
spans a 1-dimensional Toeplitz kernel if and only if $f^2$ is rigid in $H^1$ \cite{sarason94}.

The following result generalises Theorems 3.7 and 3.10 in \cite {CP}, see also \cite{Nak}.

\begin{thm}\label{2.7A}
If $g\in L^\infty(\TT)$ admits a factorisation
\begin{equation}\label{f1}
g=g_-\,\theta^{-N}g_+^{-1}
\end{equation}
where $\overline{g_-}$ and $g_+$ are outer functions in $H^2$, $g_+^2$ is rigid in $H^1$, $\theta$ is an inner function and $N\in\mathbb Z$, then
\[
\ker T_g\neq \{0\} \Leftrightarrow N>0.
\]
If $N>0$ and $\theta$ is a finite Blaschke product of degree $n$, then $\dim \ker T_g=nN$; if $\theta$ is not a finite Blaschke product, then $\dim \ker T_g=\infty$.
\end{thm}

\beginpf
(i) For $N<0$, it follows from Theorem 3.7 in \cite{CP} (proved in the context of $L^2(\RR)$) that $\ker T_g=\{0\}$.

(ii) If $N=0$, we have $g=g_-\,g_+^{-1}$ and $\ker T_g$ consists of the functions $\phi_+\in H^2$ such that $g\phi_+=\bar z\,\overline{\psi_+}$ with $\psi_+\in H^2$. We have
\begin{equation}\label{f2}
g_-\,g_+^{-1}\phi_+=\bar z\,\overline{\psi_+}\Leftrightarrow \bar z\, \frac {g_-}{\overline{g_+}} \frac{\overline{g_+}}{g_+} \phi_+=\bar z^2\,\overline{\psi_+}\Leftrightarrow \bar z  \frac{\overline{g_+}}{g_+} \phi_+=\bar z^2\frac {\overline{g_+}}{g_-}\,\overline{\psi_+}.
\end{equation}
The left-hand side of the last equality belongs to $L^2(\TT)$ while the right-hand side belongs to $\bar z^2\overline{\mathcal N_+}$, so we conclude that $\bar z^2\frac {\overline{g_+}}{g_-}\,\overline{\psi_+}\in \bar z^2\overline{H^2}\subset \overline{H_0^2}$ and, therefore, $\phi_+\in \ker T_{\bar z  \frac{\overline{g_+}}{g_+}}$. Since $g_+^2$ is rigid in $H^1$, $\ker T_{\bar z  \frac{\overline{g_+}}{g_+}}= \spam \{g_+\}$ (\cite {sarason94}): thus $\phi_+=Ag_+$ with $A\in \CC$. Now from the last equality in \eqref{f2} it follows that $Ag_+=\bar z \overline{\psi_+}$, so we cannot have $\overline{g_-}$ outer in $H^2$ unless $A=0$, i.e., $\phi_+=0$.

(iii) let now $N>0$. We have
\[
g\phi_+\in \overline{H_0^2} \Leftrightarrow   g_-\,\theta^{-N}g_+^{-1} \phi_+ \in \overline{H_0^2};
\]
any function $\phi_+= g_+\, k_a^\theta$, with $|a|<1$, satisfies that condition and therefore belongs to $\ker T_g$. This shows that $\ker T_g\neq \{0\}$ and $\dim \ker T_g=\infty$ if $\theta$ is not a finite Blaschke product. If $\theta$ is a finite Blaschke product of degree $n$, then $\theta=h_- \,z^n h_+$ with rational  left and right factors $h_\pm\in \mathcal G H^\infty$; it then follows from Theorem 3.7 in \cite{CP} that $\dim \ker T_g=nN$.
\endpf

\begin{ex}
{\rm Let $g=\frac{(z-1)^{8/15}}{z^2}\,\,,\,\,h=\frac{(z-1)^2(z+1)^{1/5}}{z^4}$ where the branches of $(z-1)^{8/15}$ and $(z+1)^{1/5}$ are analytic in $\DD$. We have 
\[
\ker T_g= \spam \{(z-1)^{7/15}\}\,\,,\,\, \ker T_h= \spam \{(z+1)^{4/5}\,,\,(z+1)^{-1/5}\}
\]
 and
\[
\bar z g^{-1} h=g_-\bar \theta g_+^{-1}\,.
\]
where $g_-=1-\bar z$ is such that $\overline{g_-}\in H^2$ is outer, $g_+ =\frac{(z-1)^{8/15}}{(z-1)(z+1)^{1/5}}\in H^2$ is such that $g_+^2$ is rigid (because $\ker T_{\bar z  \frac{\overline{g_+}}{g_+}}= \spam \{g_+\}$) and $\theta=z^2$. By solving the Riemann-Hilbert problem
\[
\bar z\,g^{-1} h\,\phi_+=\bar z \overline{\psi_+}
\]
with $\psi_+\in H^2$, we obtain
\[
\ker T_{\bar z\,g^{-1} h}=\left \{\frac{Az+B}{(z-1)^{7/15}(z+1)^{1/5}}\,:\,A, B \in \CC\right\}
\]
\[
=\spam \left\{\frac{(z-1)^{8/15}}{(z+1)^{1/5}}\,,\,\frac{1}{(z-1)^{7/15}(z+1)^{1/5}}\right\}.
\]
From Corollary \ref{2.3A} it follows that
\[
\mathcal M_2 (\ker T_g\,,\,\ker T_h)= \spam \left\{\frac{(z-1)^{8/15}}{(z+1)^{1/5}}\right \}.
\]
}
\end{ex}

The representation \eqref {f1} generalises the so called $L^2$- factorisation, which is a representation of $g$ as a product
\begin{equation}\label{f3}
g=g_-\,d\,g_+^{-1}
\end{equation}
where $g_+^{\pm1} \in H^2\,,\,g_-^{\pm1} \in \overline{H^2}$ and $d=z^k\,,\,k\in \mathbb Z$ (\cite{LS}. If $g$ is  invertible in $L^\infty(\TT)$ and admits an $L^2$-factorisation, then $\dim \ker T_g=|k|$ if $k\leq 0$, $\dim \ker T_g^*=k$ if $k\geq 0$. The factorisation \eqref{f3} is called a bounded factorisation when ${g_+}^{\pm 1} \,,\,\overline{g_-^{\pm 1}} \in {H^\infty}$. 
In various subalgebras of $L^\infty(\TT)$, every invertible element admits a factorisation \eqref{f3} where $d$ is an inner function (\cite{LS}). This is the case of the algebra of functions continuous on $\TT$ (including all rational functions without zeroes or poles on $\TT$) and the algebra $AP$ of almost periodic functions on the real line. 
In the latter case $d$ is a singular inner function, $d(\xi)=\exp (-i\lambda\xi)$ with $\lambda \in \RR$ (\cite{CD},\cite{FG}), and we have that if $g\in AP$ is invertible in $L^\infty (\RR)$ then $\ker T_g$ is either trivial or isomorphic to an infinite dimensional model space $K_\theta$ with $\theta(\xi)=\exp (i\lambda \xi)$, depending on whether $\lambda\leq 0$ or $\lambda>0$.

Various results regarding the dimension of $\ker T_{\bar z\theta h}$ can also be found in \cite{CP}  and \cite{CMP}. Namely,  if $\theta$ is a finite Blaschke product,   $\ker T_{\bar z\theta h}$ and $\ker T_{\bar z h}$ are both finite dimensional or not and, for $\dim \ker T_{\bar z h}<\infty$, we have
\[
\dim \ker T_{\bar z\theta h}=\max \{0, \dim \ker T_{\bar z h}-k\},
\]
where $k$ is the degree of $\theta$ (\cite{CMP} Theorem 6.2).

\begin{ex}{\rm For $\theta (z)=\exp(\frac{z+1}{z-1})\,,\,\phi(z)=\exp(\frac{z-1}{z+1})$, we have $\ker T_{\bar z\theta \bar \phi}=\{0\}$ (\cite{CMP}, Example 6.3); therefore $\M(K_\theta,K_\phi)=\{0\}$.}
\end{ex}

For two inner functions $\phi, \theta \in H^\infty$ we write $\phi \preceq \theta$ if $\phi$ divides $\theta$ in $H^\infty$;
that is, $\theta=\phi\psi$ for some $\psi \in H^\infty$. If we have strict inequality, that is, $\phi$ divides $\theta$
but not conversely, then we write $\phi \prec \theta$.

\begin{ex}{\rm Let $\theta\,,\,\phi$ be two inner functions with $\phi\preceq\theta$ (the case $\theta\prec\phi$ will be considered in Example \ref{2.9a}). Then $\dim \ker T_{\bar z\theta \bar \phi} \leq 1$, since $\theta \bar \phi \in H^\infty$ and $\ker T_{\theta \bar \phi}=0$  (see \cite{BCD}). We have $ \ker T_{\bar z\theta \bar \phi} =\CC$ if $\phi=a\theta$ with $a\in\CC\,,\,|a|=1$, and $ \ker T_{\bar z\theta \bar \phi} =\{0\}$ if $\phi\prec\theta$. Therefore $\M(K_\theta,K_\phi)\neq\{0\}$ if and only if $K_\theta=K_\phi$, in which case $\M(K_\theta,K_\phi)=\CC$. } 
\end{ex}

In \cite{FHR} there is a supplementary theorem describing $\M_\infty(K_\theta,K_\phi) = \M(K_\theta,K_\phi) \cap H_\infty$.
Starting with Theorem \ref{thm:main},
we immediately  have the following general result on
noting that the Carleson
measure condition is   redundant for bounded $w$.

\begin{cor}\label{bd}
Let $g,h \in L^\infty(\TT)\setminus\{0\}$ such that $\ker T_g$ and $\ker T_h$ are nontrivial. Then 
the following conditions are equivalent.\\
(i) $w \in \M_\infty(\ker T_g,\ker T_h)=\M(\ker T_g,\ker T_h) \cap H^\infty$;\\ 
(ii) $w \in H^\infty$ and $w k \in \ker T_h$ for some maximal vector $k \in \ker T_g$;\\
(iii)  $w \in H^\infty$ and $whg^{-1}\in \overline{H^\infty}$ (assuming $hg^{-1}\in L^\infty(\TT)$).\\
If $w \in H^2$,
\[
w \in \M_\infty(\ker T_g,\ker T_h)\Leftrightarrow w\in \ker T_{\bar z hg^{-1}}\cap H^\infty
\]
and if moreover $\ker T_g $ contains a maximal vector  $k$ with $k,k^{-1} \in L^\infty(\TT)$, then
\[
w \in \M_\infty(\ker T_g,\ker T_h)\Leftrightarrow wk \in \ker T_h \cap H^\infty.
\]
\end{cor}

For model spaces, we therefore recover the main theorem on bounded multipliers from \cite{FHR}.

\begin{cor}\label{fhr} \cite{FHR}
Let $\theta$ and $\phi$ be inner functions and let $w \in H^2$. Then the following
are equivalent:\\
(i) $w \in \M_\infty (K_\theta,K_\phi)$;\\
(ii) $w \in \ker T_{\overline\phi \theta \overline z} \cap H^\infty$;\\
(iii) $w S^*\theta \in K_\phi \cap H^\infty$;\\
(iv) $w\in H^\infty$ and $\bar\phi\,\theta\,w\in \overline{H^\infty}$.
\end{cor}

\beginpf
The equivalence of (i) and (ii) is contained in Corollary \ref{2.3A}. The equivalence with
(iii) follows since $S^*\theta$ is a maximal vector for $K_\theta$ that is invertible in $L^\infty(\TT)$ and the equivalence with (iv) follows from Corollary \ref{bd} (iii).
\endpf

\begin{ex}\label{2.9a}
{\rm  Let $\theta \prec \phi$; then $\ker T_{\bar z\theta \bar \phi}=K_{ z\bar\theta \phi}$ and we have $\M_\infty(K_\theta,K_\phi)=K_{ z\bar\theta \phi}\cap H^\infty$. If $\phi$ is a finite Blaschke product, then 
\[
\M_2(K_\theta,K_\phi)=\M_\infty(K_\theta,K_\phi)=K_{ z\bar\theta \phi}.
\]
}
\end{ex}

\begin{ex}
{\rm It is easy to see that a function $w_+\in H^\infty$, with an inverse in the same space, is a bounded multiplier for Toeplitz kernels. Namely, $w_+\,\ker T_g=\ker T_{g\,w_+^{-1}} \subset \ker T_{g\,w_+^{-1}\,f_-}$ for any $g\in L^\infty (\TT)\,,\,f_-\in \overline{H^\infty}$. }
\end{ex}

Applying the results of Corollary \ref{bd} to $w=1$, we also have:

\begin{prop} \label{inclusion}
Let $g,h \in L^\infty(\TT)\setminus\{0\}$, such that $\ker T_g$ and $\ker T_h$ are nontrivial. Then the following conditions are equivalent.\\
(i) $\ker T_g\subset \ker T_h$;\\
(ii) $hg^{-1}\in \overline{\mathcal N_+}$;\\
(iii) there exists a maximal function for $\ker T_g$, $k$, such that $k\in \ker T_h$.

If moreover $\ker T_g $ contains a maximal vector  $k$ with $k,k^{-1} \in L^\infty(\TT)$, then each of the above conditions is equivalent to\\
(iv) $k\in \ker T_h\cap H^\infty$. 
\end{prop}

\begin{cor}\label{2.12}
With the same assumptions as in Proposition \ref{inclusion}, if \, $hg^{-1}\in L^\infty(\TT)$, then 
\[
 \ker T_g\subset \ker T_h \Leftrightarrow  hg^{-1}\in \overline{H^\infty}
\]
\end{cor}

\begin{rem}
{\rm
Assuming without loss of generality that $hg^{-1}\in L^\infty(\TT)$, we see from the corollary above that if $\ker T_g\subset \ker T_h$ then $h=g \,  \overline{f_+}$ with $f_+ \in H^\infty$.  Let $\theta$ denote the inner factor of $f_+$. Since $\ker T_h=\ker T_{g  \overline{f_+}}=\ker T_{g \,\bar \theta}$, denoting $g \bar \theta=\tilde g$ we conclude that a Toeplitz kernel is contained in another Toeplitz kernel if and only they take the form $\ker T_{\tilde g}$ and $\ker T_{\theta\,\tilde g}$ respectively, for some inner $\theta$ and $\tilde g \in L^\infty(\TT)$.}
\end{rem}

\begin{cor}\label{2.12}
Let $g,h \in L^\infty(\TT)\setminus\{0\}$, such that $\ker T_g$ and $\ker T_h$ are nontrivial. Then $\ker T_g=\ker T_h$
if and only if there are outer functions $p,q \in H^2$ such that $\ds \frac{g}{h}=\frac{\overline p}{\overline q}$.\\
If moreover $hg^{-1}\in \mathcal G L^\infty(\TT)$, we have
\[
 \ker T_g= \ker T_h \Leftrightarrow \overline{ hg^{-1}}\in \mathcal G H^\infty.
\]
\end{cor}

It follows from Corollary \ref{2.12}, in particular, that if $h\in L^\infty(\TT)$ then $\ker T_h$ is a model space $K_\theta$ if and only if $h=\theta h_-$ with $h_-\in \mathcal G \overline{H^\infty}$.\\

In view of Corollary \ref{2.12}, one may also ask which Toeplitz kernels are contained in a model space and vice-versa.

Regarding the first question, it is clear that if $g\in \mathcal G L_\infty(\TT)$ and $\theta$ is an inner function, then $\ker T_g\subset K_\theta$ if and only if 
\beq\label {A}
g=\overline{ \theta({f_+}^{-1})}\quad {\rm  with}\quad f_+ \in H^\infty.
\eeq
If $f_+=\alpha O$ is an inner-outer factorisation with $\alpha$ inner and $O$ an outer function, from \eqref{A} we see that $\bar O \in \mathcal G \overline{H^\infty}$ because ${\bar O}^{-1}=g \theta \bar \alpha \in \overline{\mathcal N_+}\cap L^\infty(\TT)=\overline{H^\infty}$ and therefore we must have $\ker T_g=\ker T_{\bar \theta \alpha}$. In particular if $g=\bar \alpha$ where $\alpha$ is an inner function, we get the known relation  $K_\alpha\subset K_\theta \Leftrightarrow  \alpha \preceq \theta$.

Regarding the second question, we have $ K_\theta\subset\ker T_g$ with $g\in L^\infty(\TT)$ if and only if $g \in \overline{\theta H^\infty}$. In particular if $g=\bar \phi$ where $\phi$ is an inner function, we get the known relation  $K_\theta\subset K_\phi \Leftrightarrow \theta \preceq \phi$.\\

\begin{ex}{\rm
Let $\theta(z)=z^2$, so that $K_\theta=\ker T_{\bar z^2}$ is the 2-dimensional space spanned by $1$ and $z$.
The maximal vectors for this Toeplitz kernel have the form $k=a+bz$, where $ \theta \overline z \overline{a+bz}$ is outer.
That is, $\overline a z+ b $ is outer, so $0 \le |a| \le |b|$ (we should exclude the case $a=b=0$).

In other words, the non-trivial Toeplitz kernels properly contained in $K_\theta$ are 1-dimensional
and spanned by functions  $1+bz$ with $|b| < 1$, of the form $(1+bz)K_z=\ker T_{(\bar z)^2 \frac{z+\bar b}{1+bz}}$ where $\frac{z+\bar b}{1+bz}$ is an inner function.
For $b=0$ we obtain the model space $K_z$.

Note that for the non-maximal vectors $f(z)=1+bz$ for $|b|<1$ the function $w(z)=1/(1+bz)$ 
satisfies $wf \in K_\theta$, 
and $|w|^2 \, dm$ is a Carleson measure for $K_\theta$; however
 $w$ does not multiply $K_\theta$ into itself.  
}
\end{ex}

Using Proposition \ref{inclusion} and the previous results, we can study in particular the multipliers for Toeplitz kernels related by inclusion.

\begin{prop}\label{2.13}
Let $g,h \in L^\infty(\TT)\setminus\{0\}$, with $hg^{-1}\in L^\infty(\TT)$.\\
(i) If $\ker T_g\subset \ker T_h$, then
\[
\M_2 (\ker T_g,\ker T_h) = \mathcal C(\ker T_g)\cap K_{z\alpha}
\]
where $\alpha$ is the inner factor in an inner-outer factorisation of $\overline{hg^{-1}}\in H^\infty$.\\
(ii) If $\ker T_h\subset \ker T_g$, then $\M_2 (\ker T_g,\ker T_h)=\{0\}$ unless $\ker T_g= \ker T_h$.  
\end{prop}

\beginpf
(i) If $\ker T_g\subset \ker T_h$ then, by Corollary \ref{inclusion}, $hg^{-1}=\overline{f_+} \in \overline{H^\infty}$. Let $\alpha$ and $O$ denote the inner and outer factors of $f_+$, respectively. Since $\ker T_{\bar z \,\overline{f_+}}=\ker T_{\bar z \bar \alpha}$, we have from Corollary \ref{2.3A} that
\[
w\in \M_2 (\ker T_g,T_h) \Leftrightarrow w \in  \mathcal C(\ker T_g) \cap K_{z\alpha}.
\]
(ii) If $\ker T_h\subset \ker T_g$, then $hg^{-1}=(\overline{f_+})^{-1}$ with $ f_+\in H^\infty$. We have
\[
w\in \ker T_{\bar z\,(\overline{f_+})^{-1}}  \Leftrightarrow w\in H^2\,,\,\bar z\,(\overline{f_+})^{-1} w=f_-\in \overline{H^2_0}.
\]
Since $f_-\,\overline{f_+}\in  \overline{H^2_0}$, it follows that $\bar z w\in  \overline{H^2_0}$, i.e. $w\in K_z=\CC$. If $w=A\in \CC \setminus \{0\}$, then $f_+ \in \CC \setminus \{0\}$ because
\[
\bar z A=f_-\,\overline{f_+} \Rightarrow A=\overline{f_+}(zf_-)\;\;{\rm with}\,\, zf_-\in \overline{H^2}
\]
and, from the uniqueness of the inner-outer factorisation (modulo constants) it follows that $f_+$ is a constant.
\endpf

\begin{ex}
{\rm Let $\alpha$ and $\theta$ be inner with $\alpha \prec \theta$; then $\M_2 ( K_\theta,K_\alpha)=\{0\}$ and $\M_2 ( K_\alpha\,,\,K_\theta)=\mathcal C(K_\alpha)\cap K_{z\,\theta\,\bar\alpha}$. For instance, if $\theta=z^m\,,\,\alpha=z^n$ with $n\leq m$, then $\M ( K_{z^n},\,K_{z^m})=\M_2 ( K_{z^n},\,K_{z^m})=\M_\infty ( K_{z^n},\,K_{z^m})= K_{z^{m-n+1}}$}.
\end{ex}

We can generalise the results of Propositions \ref{inclusion} and \ref {2.13} for Toeplitz kernels that are equivalent in a certain sense (\cite {CMP}).

\begin{defn}\label{g}
If $g_1\,,\,g_2 \in L^\infty(\TT)$, we say that $g_1 \sim g_2$ if and only if there are functions $h_+ \in \mathcal{G}H_\infty\,,\,h_- \in \mathcal{G}\overline{H_\infty}$,  such that
\begin{equation}\label{3.22}
g_1=h_-g_2h_+.
\end{equation}
\end{defn}

It is easy to see that we have $g_1=h_-g_2h_+$ and $g_1=\tilde{h}_-g_2\tilde{h}_+$ with $h_+ \,,{\tilde h}_+\in \mathcal{G}H_\infty$ and $h_- \,,{\tilde h}_-\in \mathcal{G}\overline{H_\infty}$,  if and only if $\frac{h_-}{\tilde{h}_-}=\frac{\tilde{h}_+}{h_+}=c \in \mathbb{C}\setminus \{0\}$. If $|g_1|=|g_2|=1$ we can choose $h_{\pm}$ in \eqref{3.22} such that $\|h_-\|_\infty=\|h_+\|_\infty=1$.



\begin{defn}\label{def:3.11}
If $g_1\,,\,g_2\in L^\infty(\TT)\setminus\{0\}$, such that 
$\ker T_{g_1}\,,\,\ker T_{g_2}$ are nontrivial, we say that $\ker T_{g_1}\sim \ker T_{g_2}$ if and only if
\begin{equation}\label{3.25}
\ker T_{g_1}=h_+\ker T_{g_2}\ \quad \text{with } h_+ \in \mathcal{G}H^\infty.
\end{equation}
\end{defn}

It is clear that
$g_1 \sim g_2 \Rightarrow \ker T_{g_1}\sim \ker T_{g_2}$
since
\[
\ker T_{g_1}=\ker T_{h_-g_2h_+}=h_+^{-1}\ker T_{g_2}.
\]
It follows from Corollary \ref{2.12} that, if $g_1{g_2}^{-1} \in\mathcal G L^\infty(\TT)$, the converse is true since
\[
\ker T_{g_1}=h_+^{-1}\ker T_{g_2}\Leftrightarrow \ker T_{g_1}=\ker T_{g_2h_+}\Leftrightarrow g_1\,{g_2}^{-1}{h_+}^{-1}\in \mathcal G\overline{ H^\infty}.
\]
Therefore, if $h_+\in\mathcal G H^\infty$,
\begin{equation}\label{3.26}
\ker T_{g_1}=h_+^{-1}\ker T_{g_2}\Leftrightarrow g_1=h_-\,g_2\,h_+\quad\quad {\rm with}\;h_-\in\mathcal G \overline{H^\infty.}
\end{equation}

If $\theta_1$ is a finite Blaschke product, then it is easy to see that  $\theta_1=h_-\,z^{N_1}h_+$ where $h_+ \in \mathcal{G}H_\infty\,,\,h_- \in \mathcal{G}\overline{H_\infty}$ are rational and $N_1$ is the degree of $\theta_1$. Thus  $\theta_1\sim z^{-N_1}$. We have $K_{\theta_1 }\sim K_{ \theta_2}$ if and only if $\theta_2$ is also a finite Blaschke product of the same degree. Moreover, if $\theta_1$ and $\theta_2$ are finite Blaschke products with $\theta_1\sim z^{-N_1}$ and $\theta_2\sim z^{-N_2}$, then $\overline {\theta_1}\,\theta_2 \sim z^{N_1-N_2}$ and we have 
\[
\ker T_{\overline {\theta_1}\,\theta_2 }=\{0\}\quad {\rm if} \; N_2\leq N_1\;,\;\ker T_{\overline {\theta_1}\,\theta_2 }\sim K_{z^{N_1-N_2}} \quad {\rm if} \; N_1> N_2.
\]

\begin{prop}
Let $g,h \in L^\infty(\TT)\setminus\{0\}$, with $hg^{-1}\in L^\infty(\TT)$.\\
(i) $\ker T_g\sim\ker T_{\tilde g}\subset \ker T_h$ for some $\tilde g \in L^\infty(\TT)
$ if and only if there exists $h_+\in \mathcal G H^\infty$ such that $hg^{-1}h_+ \in \overline{H^\infty}$.\\
(ii) If  $\ker T_g\sim\ker T_{\tilde g}\subset \ker T_h$ for some $\tilde g \in L^\infty(\TT)$, with $\ker T_g= {h_+}^{-1}\ker T_{\tilde g}$ where $h_+ \in \mathcal G H^\infty$, then
\[
\M_2(\ker T_g, T_h)=h_+^{-1}\M_2(\ker T_{\tilde g},\ker T_h)=\mathcal C(\ker T_g)\cap h_+ K_{z\alpha}
\]
where $\alpha$ is the inner factor of an inner-outer factorisation of $\overline{hg^{-1}h_+}\in H^\infty$.
\end{prop}

\beginpf
(i) If $\ker T_g \sim \ker T_{\tilde g}$ then by Definition \ref{def:3.11} and \eqref{3.26} there exist $h_+ \in \mathcal{G}H_\infty\,,\,h_- \in \mathcal{G}\overline{H_\infty}$,  such that
$g=h_- \tilde g h_+$; on the other hand, by Corollary \ref{2.12}
\[
 \ker T_{\tilde g}\subset \ker T_h \Leftrightarrow  h{\tilde g}^{-1}\in \overline{H^\infty}\Leftrightarrow  hh_-{ g}^{-1}h_+\in \overline{H^\infty}\Leftrightarrow  h{ g}^{-1}h_+\in \overline{H^\infty}.
\]
Conversely, if there exists $h_+\in \mathcal G H^\infty$ such that $h{ g}^{-1}h_+\in \overline{H^\infty}$, then $ \ker T_{ gh_+^{-1}}\subset \ker T_h$ and taking $\tilde g=gh_+^{-1}$ we conclude that $\ker T_g\sim \ker T_{\tilde g}\subset \ker T_h$.

(ii) If $\ker T_g= {h_+}^{-1}\ker T_{\tilde g}$, we have $\M(\ker T_g,\ker T_h)={h_+}^{-1} \M(\ker T_{\tilde g,\ker T_h})$ and by Proposition \ref{2.13}
\[
\M_2 (\ker T_{\tilde g},\ker T_h) = \mathcal C(\ker T_{\tilde g})\cap K_{z\alpha}
\]
where $\alpha$ is the inner factor of $\overline{h{\tilde g}^{-1}}\in H^\infty$, which is equal to the inner factor of $\overline{hg^{-1}h_+}\in H^\infty$.
\endpf

\section{Surjective multipliers}
\label{sec:3}

The original context of Crofoot's work \cite{crofoot}
is where the multiplication operator between two model spaces is surjective. We may 
obtain similar results in the more general context of Toeplitz kernels.

\begin{lem}\label{lem:feb23}
Let $g \in L^\infty(\TT)$, 
let $k$ be a maximal vector for $\ker T_g$, and suppose that  $w \ker T_g$ is a Toeplitz kernel. Then   $w \ker T_g=\Kmin(wk)$.
\end{lem}

\beginpf
Let $h \in L^\infty(\TT)$ be such that $w \ker T_g=\ker T_h$.
We have $wk \in \ker T_h$ and $\ker T_h=w \ker T_g \subset \Kmin(wk)$ by Theorem \ref{thm:main}. 
Hence  $\ker T_h=\Kmin(wk)$.
\endpf

\begin{thm}\label{thm:3.1}
Let $g, h \in L^\infty(\TT)$ such that $\ker T_g$ and $\ker T_h$ are nontrivial.  
Then a function $w \in \Hol(\DD)$ satisfies $w \ker T_g=\ker T_h$
if and only if\\
(i) $w \in \mathcal C (\ker T_g)$ and $w^{-1}\in \mathcal C ( \ker T_h)$;\\
(ii) for some (or indeed, for every) maximal vector $k \in \ker T_g$, the function $wk$ is a maximal
vector for $\ker T_h$.
\end{thm}

\beginpf
Suppose that the conditions are satisfied. Then by Theorem \ref{thm:main} $w$ is a multiplier from
$\ker T_g$ into $\ker T_h$ and $w^{-1}$ is a multiplier from $\ker T_h$ into $\ker T_g$. Since
the multiplication operator is injective, we see that we have $w \ker T_g = \ker T_h$.

Conversely, if $w \ker T_g=\ker T_h$, then condition (i) is clearly satisfied, and (ii) follows from Lemma~\ref{lem:feb23}.
\endpf





We also have the following necessary and sufficient condition:

\begin{thm}\label{thm:3.4}
Let $g, h \in L^\infty(\TT)$ such that $\ker T_g$ and $\ker T_h$ are nontrivial. 
Then $w \ker T_g=\ker T_h$ if and only if 
$w \in \mathcal C (\ker T_g)\,,\,w^{-1}\in \mathcal C ( \ker T_h)$ and
\beq\label{eq:hlinkedtog}
h= g \frac{\overline w}{w} \frac{\overline q}{\overline p}
\eeq
for some outer functions $p, q \in H^2$.
\end{thm}
\beginpf
Note that $w$ must be outer, as functions in a Toeplitz kernel cannot share a common
inner factor, since if $f \in \ker T_g$  and $\theta$ is inner with $f/\theta \in H^2$, then $f/\theta \in \ker T_g$.

Now let $k=\theta u$ be a maximal vector for $\ker T_g$, where $\theta$ is inner and $u$ is outer.
Then $\ker T_g = \ker T_{\overline z\overline \theta \overline u/u}$. 
We write $g_0=\overline z\overline \theta \overline u/u$.
Also   the inner--outer factorization of $wk$, which is a maximal vector for $\ker T_h$,
is
$wk=\theta (w u)$, so we have $\ker T_h=\ker T_{\overline z\overline \theta \overline w\overline u/(wu)}$.
We write $h_0=\overline z\overline \theta \overline w\overline u/(wu)$.

By Corollary \ref{2.12} we have outer functions $r$ and $s$ such that
$g= g_0 \overline r/\overline s$.  So
\[
\ker T_h =  \ker T_{h_0} =\ker T_{g_0 \overline w/w} = \ker T_{g_0 \overline w \overline r/(w\overline s)} = \ker T_{g \overline w/w}.
\]
Finally, by Corollary \ref{2.12} we have \eqref{eq:hlinkedtog}.

For the converse, we see that \eqref{eq:hlinkedtog} implies that $\ker T_h = \ker T_{g\overline w/w}$.
Then if $f \in \ker T_g$ we have $(fw)(g \overline w/w)= fg \overline w \in \H2p$ and so $fw \in \ker T_{g\overline w/w}=\ker T_h$.
Also if $f \in \ker T_h$ then $fg/w = (fg\overline w/w)/\overline w \in \H2p$, and so $f/w \in \ker T_g$.
\endpf

\begin{rem}{\rm
In the case of model spaces, suppose that $wK_\theta=K_\phi$;
then we apply the above results to $g=\overline\theta$ and $h=\overline\phi$, so
we have $K_\phi = \ker T_{\overline \theta \overline w/w}$.
Now $\theta     w/\overline w \in L^\infty(\TT)$ (indeed it is unimodular), but it
also equals $\phi p/q$ from \eqref{eq:hlinkedtog}, and this is in the Smirnov class;
so it lies in $H^\infty$ and is inner.

Thus $K_\phi=K_{\theta w/\overline w}$, and so $\phi= \alpha \theta w/\overline w$, with $\alpha \in \CC$ 
and $|\alpha|=1$,
which is
Crofoot's result.}
\end{rem}

The equivalence relation of Definition \ref{def:3.11} is closely related to the question of existence of surjective multipliers between two Toeplitz kernels. Indeed, any $w=w_+\in \mathcal G H^\infty$ is a surjective multiplier from any given $\ker T_g$ onto another Toeplitz kernel $\ker T_{{w_+}^{-1}g} = w_+ \ker T_g$. One may ask if the same is true for model spaces, i.e., given $w_+ \in \mathcal G H^\infty$ and an inner function $\theta$, is there always another inner function $\phi$ such that $w_+ K_\theta\subset K_\phi$?

The answer to this question is negative. In fact, if $\theta$ is a finite Blaschke product then $K_\theta=\ker T_{\bar \theta}$ and $w_+ K_\theta=\ker T_{{w_+}^{-1}\bar\theta}$ must both be finite dimensional, with the same dimension. If $w_+ K_\theta=K_\phi$ with $\phi$ inner, then  we must have, on the one hand, $w_+\theta\bar\phi\in \mathcal G \overline{H^\infty}$ and on the other hand, since $\theta\sim z^{-N}\,,\, \phi\sim z^{-N}$ for some $N\in \mathbb N$, we must have $h_-w_+h_+ =f_-$ for some rational $h_-\in \mathcal G \overline{H^\infty}\,,\, h_+\in \mathcal G H^\infty$ and $f_-\in \mathcal G \overline{H^\infty}$. It follows that $w_+h_+ =A \in \CC$ and therefore $w_+ K_\theta=K_\phi$ only if $w_+$ is a rational function in $\mathcal G H^\infty$.


\section{The upper half-plane}
\label{sec:4}

The results on Toeplitz kernels in \cite{CP14,CMP} were originally derived for the
Hardy space $H^2(\CC^+)$ of the upper half-plane.
There are additional motivations
here, in that Paley--Wiener spaces appear naturally in the context of model
spaces corresponding to the inner functions $\theta(s)=e^{i \lambda s}$ for $\lambda > 0$:
for this and other motivations  we refer to the introduction of \cite{CP17}.

Recall that we have the relation $H^2(\CC^-)=L^2(\RR) \ominus H^2(\CC^+)$, and
$f \in H^2(\CC^-)$ if and only if $\overline f \in H^2(\CC^+)$.

Moreover it is well known (see, e.g. \cite[pp. 23--24]{Par}) that $g \in L^p(\RR)$ 
for some $1 \le p < \infty$ if and only if the function
$V_p g$ defined by 
\beq\label{eq:isom}
V_pg(z)=2^{2/p} \pi^{1/p} (1+z)^{-2/p} g(i(1-z)/(1+z))
\eeq
lies in $L^p(\TT)$. Indeed, $V_p$ is an isometric map which preserves the
corresponding Hardy spaces,  with $H^p(\CC^+)$ mapping to $H^p(\DD)$.

The analogue of Theorem \ref{thm:main} is the following. We now use $m$ to
refer to Lebesgue measure on $\RR$, and $T_g$ etc. to refer to Toeplitz operators
on $H^2(\CC^+)$.

\begin{thm}\label{thm:main2}
Let $g, h \in L^\infty(\RR)$ such that $\ker T_g$ and $\ker T_h$ are nontrivial.  
Then a function $w \in \Hol(\CC^+)$ lies in $ \M(\ker T_g, \ker T_h)$ if and only if\\
(i) $w k \in \ker T_h$ for some (and hence all) maximal vectors $k$ of $\ker T_g$;\\
(ii) $w \ker T_g \subset L^2(\RR)$; that is $|w|^2 \, dm$ is a Carleson measure for $\ker T_g$.
\end{thm}
\beginpf
Clearly, the two conditions are necessary. So assume that (i) and (ii) hold, and write $k=\theta p$, where
$\theta$ is inner and $p$ is outer.
Now $\ker T_g= \ker T_{\overline \theta\overline p/p}$, as detailed above,
and
thus without loss of generality we may take $g=\overline \theta\overline p/p$.

We have that $wkh \in \Hlp$, since $wk \in \ker T_h$. Suppose now that
$f\in \ker T_g$, so that $ f g \in\Hlp$.
Now 
\[
wf  h = (wkh) \frac{f}{\theta p} = (wkh) \frac{fg}{\overline p}.
\]
Then $w f h \in L^2(\RR)$, since $wf \in L^2(\RR)$ by the Carleson condition.

Also $wkh$ and $f g$ are in $\Hlp$ so  
$ \overline{wfh}=  \overline{wkh}\, \overline{fg}/p$
is in the Smirnov class
of the half-plane (the ratio of an $H^1(\CC^+)$ function and an outer $H^2$ function) as well as $L^2(\RR)$.
The generalized maximum principle  
applies also to the half-plane, as can be seen using the isometric equivalences in \eqref{eq:isom}.
We conclude that
$ \overline{wfh} \in H^2(\CC^+)$ and so $wf h \in \Hlp$,
and finally $wf \in \ker T_h$.
\endpf

The method of proof of Theorem \ref{thm:maxtg} shows that the maximal vectors for a nontrivial Toeplitz
kernel $\ker T_g \subset H^2(\CC^+)$ are functions of the form $g^{-1}\overline p$, where $p \in H^2(\CC^+)$ outer.
Maximal vectors for model spaces $K_\theta= \ker T_{\overline \theta}$ have already been characterized in \cite[Thm 5.2]{CMP}
as functions in $H^2(\CC^+)$ of the form $\theta \overline p$ with $p$ outer.
One such is $k(s)=(\theta(s)-\theta(i)) / (s-i)$,
the backward shift of the function $\theta$, although $\theta$ itself is not in $H^2(\CC^+)$. Since $k(s) =\theta(s) (1-\theta(i)\overline {\theta(s)})/(s-i)$ 
for $s \in \RR$ we see that this $k$ is an appropriate test function to use.\\

One special case of interest is when $\ker T_g$ consists entirely of bounded functions, since then
any $H^2$ function $w$ automatically satisfies the Carleson condition in Theorems
\ref{thm:main} and \ref{thm:main2}: this property is
discussed for model spaces in \cite{CMP}. For the disc, $K_\theta \subset H^\infty$
if and only if $K_\theta$ is  finite-dimensional, that is,   $\theta$ is
rational, but for the half-plane there are other possibilities, for example $\theta(s)=e^{i\lambda s}$
with $\lambda > 0$. We refer to \cite{CMP} for further details.

Finally, we remark that Theorems \ref{thm:3.1} and \ref{thm:3.4} hold in the case of the half-plane
with obvious modifications.

\subsection*{Acknowledgements}

The authors are grateful to the referee for a careful reading of the manuscript and some  
helpful comments.
This work was partially supported by\\ FCT/Portugal through  UID/MAT/04459/2013.

\end{document}